\documentclass[10pt,twoside]{article}
\usepackage{amscd,amssymb,amsxtra}
\usepackage[mathscr]{eucal}
\usepackage[all]{xy}
\usepackage{Latex-document}
\def\volumeno{I}
\CompileMatrices

\title{\bf Algebraic Topology and Modular Forms\vskip 6mm}

\author{M.~J.~Hopkins\thanks{Department of Mathematics, Massachusetts
Institute of Technology, Cambridge, MA 02139, USA. E-mail:
mjh@math.mit.edu} \vspace*{-0.5cm}}
\date{\vspace{-8mm}}

\newtheorem{thm}[equation]{Theorem}
\newtheorem{cor}[equation]{Corollary}

\newtheorem{prop}[equation]{Proposition}

\ifx\undefined\theoremstyle
    \newtheorem{definition}[equation]{Definition}
    \newenvironment{defin}{\begin{definition}}{\end{definition}}
    \newtheorem{conj}[equation]{Conjecture}
    \newtheorem{rem}[equation]{Remark}
    \newtheorem{rems}[equation]{Remarks}

\else
    \theoremstyle{definition}

    \newtheorem{defin}[equation]{Definition}

    \theoremstyle{remark}

    \newtheorem{rem}[equation]{Remark}
    
\ifx\undefined\prem \newtheorem{prem}{\bf Working Remark} \fi

\fi

\ifx\undefined\pf
\newenvironment{pf}{\bigskip{\em Proof:\/}}{\qed\medskip}
\newenvironment{pf*}[1]{\bigskip{\em #1:\/}}{\qed\medskip}
\fi

\ifx\undefined\numberwithin
\def\numberwithin#1#2{\makeatletter\@ifundefined{c@#1}{\@nocnterrr}{%
  \@ifundefined{c@#2}{\@nocnterr}{%
  \@addtoreset{#1}{#2}%
  \toks@\expandafter\expandafter\expandafter{\csname the#1\endcsname}%
  \expandafter\xdef\csname the#1\endcsname
    {\expandafter\noexpand\csname the#2\endcsname
     .\the\toks@}}}\makeatother}\fi

\ifx\undefined\qed\newcommand{\qed}{\hfil\rule{4pt}{6pt}\bigskip}\fi

\DeclareMathOperator{\ext}{Ext}

\newcommand{\Q}{{\mathbb Q}}
\newcommand{\Z}{{\mathbb Z}}
\newcommand{\R}{{\mathbb R}}

\newcommand{\zerowidth}[1]{\hbox to 0pt{\hss$\displaystyle #1$\hss}}
\newcommand{\LL}{[\mkern-2mu[}
\newcommand{\RR}{]\mkern-2mu]}

\newcommand{\C}{\mathbb C}

\ifx\undefined\eqref\newcommand{\eqref}[1]{\rm (\ref{#1})}\fi

\newcounter{thmItem}

\newenvironment{thmList}{\begin{list}%
{\rm \roman{thmItem})}{\usecounter{thmItem}
\setlength{\labelwidth}{2em}
\setlength{\itemindent}{2em}
\setlength{\leftmargin}{0pt}
\setlength{\listparindent}{0pt}
\setlength{\parsep}{0pt}
\setlength{\partopsep}{0pt}
\setlength{\itemsep}{\medskipamount}
\setlength{\topsep}{\medskipamount}
}}{\end{list}}

\newcounter{textItem}

\newcounter{condItem}

\ifx\undefined\slot
\newcommand{\slot}{\,-\,}
\fi

\DeclareMathOperator{\tmf}{tmf}
\DeclareMathOperator{\spin}{Spin}
\DeclareMathOperator{\bspin}{{\mathit{B}Spin}}
\DeclareMathOperator{\mspin}{{\mathit{M}Spin}}
\DeclareMathOperator{\spec}{Spec}

\newcommand{\moeight}{MO\langle8\rangle}
\newcommand{\cpinfty}{\cp^{\infty}}
\newcommand{\hpinfty}{\hp^{\infty}}
\newcommand{\cp}{CP}
\newcommand{\hp}{HP}
\newcommand{\boeight}{BO\langle 8\rangle}
\newcommand{\sboeight}{\bo\langle 8\rangle}
\newcommand{\einfty}{E_{\infty}}
\newcommand{\ainfty}{A_{\infty}}
\DeclareMathOperator{\gl}{\mathit{gl_{1}}}
\DeclareMathOperator{\Gl}{\mathit{Gl_{1}}}
\DeclareMathOperator{\BGl}{\mathit{BGl_{1}}}
\newcommand{\bo}{bo}

\DeclareMathOperator{\sym}{Sym}
\newcommand{\pfn}{\mathcal P}
\DeclareMathOperator{\res}{Res}
\DeclareMathOperator{\SL}{SL}
\newcommand{\mumu}{\langle \mu,\mu\rangle}

\newcommand{\ellmu}{\langle \ell,\mu\rangle}

\newcommand{\pistk}[1]{\pi^{\text{st}}_{#1}}

\newcommand{\bkn}{\Sigma^{K(n)}}

\begin{document}

\maketitle

\thispagestyle{first} \setcounter{page}{283}



\section{Introduction}

\vskip-5mm \hspace{5mm}

The problem of describing the homotopy groups of spheres has been
fundamental to algebraic topology for around 80 years.  There were
periods when specific computations were important and periods when
the emphasis favored theory.  Many mathematical invariants have
expressions in terms of homotopy groups, and at different times
the subject has found itself located in geometric topology,
algebra, algebraic $K$-theory, and algebraic geometry, among other
areas.

There are basically two approaches to the homotopy groups of
spheres. The oldest makes direct use of geometry, and involves
studying a map $f:S^{n+k}\to S^{n}$ in terms of the inverse image
$f^{-1}(x)$ of a regular value.  The oldest invariant, the degree
of a map, is defined in this way, as was the original definition
of the Hopf invariant.  In the 1930's
Pontryagin\footnote{Lefschetz reported on this work at the 1936
ICM in Oslo.}~\cite{pontryagin38_1,pontryagin38_2} showed that the
homotopy class of a map $f$ is completely determined by the
geometry of the inverse image $f^{-1}(B_{\epsilon}(x))$ of a small
neighborhood of a regular value.  He introduced the basics of
framed cobordism and framed surgery, and identified the group
$\pi_{n+k}S^{n}$ with the cobordism group of smooth $k$-manifolds
embedded in $\R^{n+k}$ and equipped with a framing of their stable
normal bundles.

The other approach to the homotopy groups of spheres involves
comparing spheres to spaces whose homotopy groups are known.  This
method was introduced by Serre~\cite{serre52:_sur_freud,
serre52:_sur_maclan, cartan52:_espaci,cartan52:_espacii} who used
Eilenberg-MacLane spaces $K(A,n)$, characterized by the property
\[
\pi_{i}K(A,n)=
\begin{cases}
A &\quad i=n \\
0 &\quad\text{otherwise.}
\end{cases}
\]
By resolving a sphere into Eilenberg-MacLane spaces Serre was able to
compute $\pi_{k+n}S^{n}$ for all $k\le 8$.

For some questions the homotopy theoretic methods have proved more
powerful, and for others the geometric methods have.  The resolutions
that lend themselves to computation tend to use spaces having
convenient homotopy theoretic properties, but with no particularly
accessible geometric content.  On the other hand, the geometric
methods have produced important homotopy theoretic moduli spaces and
relationships between them that are difficult, if not impossible, to
see from the point of view of homotopy theory.  This metaphor is
fundamental to topology, and there is a lot of power in spaces, like
the classifying spaces for cobordism, that directly relate to both
geometry and homotopy theory.  It has consistently proved important to
understand the computational aspects of the geometric devices, and the
geometric aspects of the computational tools.

A few years ago Haynes Miller and I constructed a series of new
cohomology theories, designed to isolate certain ``sectors'' of
computation.  These were successful in resolving several open
issues in homotopy theory and in contextualizing many others.
There seemed to be something deeper going on with one of them, and
in~\cite{hopkins95:_topol_witten} a program was outlined for
constructing it as a ``homotopy theoretic'' moduli space of
elliptic curves, and relating it to the Witten genus. This program
is now complete, and we call the resulting cohomology theory
$\tmf$ (for {\em topological modular forms}).  The theory of
topological modular forms has had applications in homotopy theory,
in the theory of manifolds, in the theory of lattices and their
$\theta$-series, and most recently seems to have an interesting
connection with the theory of $p$-adic modular forms. In this note
I will explain the origins and construction of $\tmf$ and the way
some of these different applications arise.


\section{Sixteen homotopy groups}
\label{sec:sixt-homot-groups}

\vskip-5mm \hspace{5mm}

By the Freudenthal suspension theorem, the value of the homotopy
group $\pi_{n+k}S^{n}$ is independent of $n$ for $n>k+1$.  This
group is $k^{\text{th}}$ stable homotopy group of the sphere,
often written $\pistk k(S^{0})$, or even as $\pi_{k}S^{0}$ if no
confusion is likely to result.  In the table below I have listed
the values of $\pi_{n+k}S^{n}$ for $n\gg0$ and $k\le 15$.

\bigskip
\begin{center}
\begin{tabular}{|c|c|c|c|c|c|c|c|c|c|}\hline
$k$ & $0$ & $1$ & $2$ & $3$ & $4$ & $5$ & $6$ & $7$ & $8$   \\ \hline
$\pi_{n+k}S^{n}$ & $\Z$ & $\Z/2$ & $\Z/2$  & $\Z/24$ & $0$ & $0$ &
$\Z/2$ & $\Z/240$ & $\Z/2\oplus \Z/2$ \\ \hline
\end{tabular}

\vskip .2in

\begin{tabular}{|c|c|c|c|c|c|c|c|c|}\hline
 $9$ & $10$ & $11$ & $12$ & $13$ & $14$ & $15$   \\ \hline
 $(\Z/2)^{3}$ & $\Z/6$ & $\Z/504$ & $0$ & $\Z/3$ & $(\Z/2)^{2}$ &
$\Z/2\oplus \Z/480$  \\ \hline
\end{tabular}
\end{center}
\bigskip
In geometric terms, the group $\pi_{n+k}S^{n}$ is the cobordism group
of stably framed manifolds, and a homomorphism from $\pi_{n+k}S^{n}$
to an abelian group $A$ is a cobordism invariant with values in $A$.
The groups in the above table thus represent universal invariants of
framed cobordism.  Some, but not all of these invariants have
geometric interpretations.

When $k=0$, the invariant is simply the number of points of the framed
$0$-manifold.  This is the geometric description of the degree of a
map.

When $k=1$ one makes use of the fact that any closed $1$-manifold is a
disjoint union of circles.  The $\Z/2$ invariant is derived from the
fact that a framing on $S^{1}$ differs from the framing which bounds a
framing of $D^{2}$ by an element of $\pi_{1}SO(N)=\Z/2$.

There is an interesting history to the invariant in dimension $2$.
Pontryagin originally announced that the group $\pi_{n+2}S^{n}$ is
trivial.  His argument made use of the classification of Riemann
surfaces, and a new geometric technique, now known as framed surgery.
He later~\cite{pontryagin50:_homot} correctly evaluated this group,
but for his corrected argument didn't need the technique of surgery.
Surgery didn't reappear in again until around 1960, when it went
on to play a fundamental role in geometric topology.  The invariant is
based on the fact that a stable framing of a Riemann surface
$\Sigma$ determines a quadratic function $\phi:H^{1}(\Sigma;\Z/2)\to
\Z/2$ whose underlying bilinear form is the cup product.  To describe
$\phi$, note that each $1$-dimensional cohomology class $x\in
H^{1}(\Sigma)$ is Poincar\'e dual to an oriented, embedded
$1$-manifold, $C_{x}$, which inherits a framing of its stable normal
bundle from that of $\Sigma$.  The manifold $C_{x}$ defines an element
of $\pistk 1S^{0}=\Z/2$, and the value of $\phi(x)$ is taken to be
this element.  The cobordism invariant in dimension $2$ is the Arf
invariant of $\phi$.

A similar construction defines a map
\begin{equation}
\label{eq:6}
\pistk{4k+2}S^{0}\to \Z/2.
\end{equation}
In~\cite{browder69:_kervair} Browder interpreted this invariant in
homotopy theoretic terms, and showed that it can be non-zero only for
$4k+2=2^{m}-2$.  It is known to be non-zero for $\pistk2S^{0}$,
$\pistk6S^{0}$, $\pistk{14}S^{0}$, $\pistk{30}S^{0}$ and
$\pistk{62}S^{0}$.  The situation for $\pistk{2^{m}-2}S^{0}$ with
$m>6$ is unresolved, and remains an important problem in algebraic
topology.  More recently, the case $k=1$ of~\eqref{eq:6} has appeared
in $M$-theory~\cite{witten97:_five_m}.  Building on this, Singer and
I~\cite{hopkins00:_quadr_m} offer a slightly more analytic
construction of~\eqref{eq:6}, and relate it to Riemann's
$\theta$-function.


Using $K$-theory, Adams~\cite{Ad:JX4} defined surjective
homomorphisms (the $d$ and $e$-invariants)
\begin{align*}
\pi_{4n-1}S^{0} &\to \Z/d_{n}, \\
\pi_{8k}S^{0} &\to \Z/2, \\
\pi_{8k+1}S^{0} &\to \Z/2\oplus\Z/2, \\
\pi_{8k+2}S^{0} &\to \Z/2.
\end{align*}
where $d_{n}$ denotes the denominator of $B_{2n}/(4n)$.  He (and
Mahowald~\cite{Mah:J}) showed that they split the inclusion of the
image of the $J$-homomorphism, making the latter groups summands.
A geometric interpretation of these invariants appears in
Stong~\cite{Stg} using $\spin$-cobordism, and an analytic
expression for the $e$-invariant in terms of the Dirac operator
appears in the work of Atiyah, Patodi and Singer~\cite{aps2,aps3}.
The $d$-invariants in dimensions $(8k+1)$ and $(8k+2)$ are given
by the mod $2$ index of the Dirac
operator~\cite{AtiSing:Ind5,atiyah69:_index_fredh}.

This more or less accounts for the all of the invariants of framed
cobordism that can be constructed using known geometric
techniques. In every case the geometric invariants represent
important pieces of mathematics.  What remains is the following
list of homotopy theoretic invariants having no known geometric
interpretation:

\bigskip
\begin{center}
\begin{tabular}{|c|c|c|c|c|c|c|c|c|c|c|}\hline
 $\cdots$& $8$ & $9$ & $10$ & $11$ & $12$ & $13$ & $14$ & $15$ \\ \hline
 $\cdots$   & $\Z/2$ & $\Z/2$ & $\Z/3$ &  &  & $\Z/3$ & $\Z/2$ & $\Z/2$  \\ \hline
\end{tabular}

\end{center}
\bigskip

\noindent This part of homotopy theory is not particularly exotic.  In
fact it is easy to give examples of framed manifolds on which the
geometric invariants vanish, while the homotopy theoretic invariants
do not.  The Lie groups $SU(3)$, $U(3)$, $Sp(2)$, $Sp(1)\times Sp(2)$,
$G_{2}$, $U(1)\times G_{2}$ have dimensions $8$, $9$, $10$, $13$,
$14$, and $15$, respectively.  They can be made into framed manifolds
using the left invariant framing, and in each case the corresponding
invariant is non-zero.  We will see that the theory of topological
modular forms accounts for {\em all} of these invariants, and in doing
so relates them to the theory of elliptic curves and modular forms.
Moreover many new invariants are defined.


\section{Spectra and stable homotopy}
\label{sec:techn-comp}

\vskip-5mm \hspace{5mm}

In order to explain the theory of topological modular forms it is
necessary to describe the basics of stable homotopy theory.

\subsection{Spectra and generalized homology}

\vskip-5mm \hspace{5mm}

Suppose that $X$ is an $(n-1)$-connected pointed space.  By the
Freudenthal Suspension Theorem, the suspension homomorphism
\[
\pi_{n+k}(X)\to \pi_{n+k+1}\Sigma X
\]
is an isomorphism in the range $k<2n-1$.  This is the {\em stable
range} of dimensions, and in order to isolate it and study only and
{\em stable homotopy theory} one works in the category of {\em
spectra}.
\begin{defin}
{\rm
(see~\cite{LMayS,HopGoerss:Mult,elmendorf95:_moder,elmendorf97:_rings,Ad:SHGH})}
A {\em spectrum} $E$ consists of a sequence of pointed spaces
$E_{n}$, $n=0,1,2,\dots$ together with maps
\begin{equation}
s_{n}^{E}:\Sigma E_{n}\to E_{n+1}
\end{equation}
whose adjoints
\begin{equation}\label{eq:4}
t_{n}^{E}:E_{n}\to \Omega E_{n+1}
\end{equation}
are homeomorphisms.
\end{defin}

A map $E\to F$ of spectra consists of a collection of maps
\[
f_{n}:E_{n}\to F_{n}
\]
which is compatible with the structure maps $t_{n}^{E}$ and
$t_{n}^{F}$.

For a spectrum $E=\{E_{n},t_{n} \}$ the value of the group
$\pi_{n+k}E_{n}$ is independent of $n$, and is written $\pi_{n}E$.
Note that this makes sense for any $n\in\Z$.   More generally, for any pointed
space $X$, the $E$-homology and $E$-cohomology groups of $X$ are
defined as
\begin{align*}
E^{k}(X) &=[\Sigma^{n}X,E_{n+k}], \\
E_{k}(X) &=\varinjlim \pi_{n+k}E_{n}\wedge X.
\end{align*}
Any homology theory is represented by a spectrum in this way, and
any map of homology theories is represented (not necessarily
uniquely) by a map of spectra.  For example, the spectrum $HA$
with $HA_{n}$ the Eilenberg-MacLane space $K(A,n)$ represents
ordinary homology with coefficients in an abelian group $A$.

\subsection{Suspension spectra and Thom spectra}

\vskip-5mm \hspace{5mm}

In practice, spectra come about from a sequence of spaces $X_{n}$
and maps $t_{n}:\Sigma X_{n}\to X_{n+1}$.  If each of the maps
$t_{n}$ is a closed inclusion, then the collection of spaces
\[
(LX)_{k}=\varinjlim \Omega^{n}X_{n+k}
\]
forms a spectrum.  In case $X_{n}=S^{n}$, the resulting spectrum is
the {\em sphere spectrum} and denoted $S^{0}$.  By construction
\[
\pi_{k}S^{0}=\pistk kS^{0}=\pi_{n+k}S^{n}\qquad n\gg0.
\]
In case $X_{n}=\Sigma^{n}X$, the resulting spectrum is the {\em
suspension spectrum} of $X$, denoted $\Sigma^{\infty}X$ (or just $X$
when no confusion with the space $X$ is likely to occur).  Its
homotopy groups are given by
\[
\pi_{k}\Sigma^{\infty}X = \pistk kX = \pi_{n+k}\Sigma^{n}X\qquad n\gg0,
\]
and referred to as the {\em stable homotopy groups of $X$}.

Another important class of spectra are {\em Thom spectra}.  Let
$BO(n)$ denote the Grassmannian of $n$-planes in $\R^{\infty}$, and
$MO(n)$ the Thom complex of the universal $n$-plane bundle over
$BO(n)$.  The natural maps
\[
\Sigma MO(n)\to MO(n+1)
\]
lead to a spectrum $MO$, the unoriented bordism spectrum.  This
spectrum was introduced by Thom~\cite{Thom}, who identified the group
$\pi_{k}MO$ with the group of cobordism classes of $k$-dimensional
unoriented smooth manifolds.  Using the complex Grassmannian instead of
the real Grassmannian leads to the {\em complex cobordism} spectrum
$MU$.  The group $\pi_{k}MU$ can be interpreted as the group of
cobordism classes of $k$-dimensional stably almost complex
manifolds~\cite{Mil:MU}.    More generally, a Thom spectrum $X^{\zeta}$
is associated to any map
\[
\zeta:X\to BG
\]
from a space $X$ to the classifying space $BG$ for stable spherical
fibrations.

The groups $\pi_{k}MO$ and $\pi_{k}MU$ have been
computed~\cite{Thom,Mil:MU}, as many other kinds of cobordism
groups.  The spectra representing cobordism are among the few
examples that lend themselves to both homotopy theoretic and
geometric investigation.

\subsection{Algebraic structures and spectra}
\label{sec:algebr-struct-spectr}

\vskip-5mm \hspace{5mm}

The set of homotopy classes of maps between spectra is an abelian
group, and in fact the category of abelian groups makes a fairly
good guideline for contemplating the general structure of the
category of spectra.  In this analogy, spaces correspond to sets,
and spectra to abelian groups.  The smash product of pointed
spaces
\[
X\wedge Y = X \times Y /(x,\ast)\sim (\ast,y)
\]
leads to an operation $E\wedge F$ on spectra analogous to the tensor
product of abelian groups.  Using this ``tensor structure'' one can
imitate many constructions of algebra in stable homotopy theory, and
form analogues of associative algebras ($A_{\infty}$-ring spectra),
commutative algebras ($E_{\infty}$-ring spectra), modules, etc.  The
details are rather subtle, and the reader is referred
to~\cite{elmendorf97:_rings} and~\cite{smith00:_symmet} for further
discussion.

The importance of refining common algebraic structures to stable
homotopy theory has been realized by many
authors~\cite{may77:_e_e,elmendorf97:_rings,waldhausen78:_algeb_k_i,
waldhausen79:_algeb_k_ii}, and was especially advocated by Waldhausen.

The theory of topological modular forms further articulates this
analogy. It is built on the work of Quillen relating formal groups
and complex cobordism.  In~\cite{Qui:FGL}, Quillen portrayed the
complex cobordism spectrum $MU$ as the universal cohomology theory
possessing Chern classes for complex vector bundles (a {\em
complex oriented} cohomology theory).  These generalized Chern
classes satisfy a Cartan formula expressing the Chern classes of a
Whitney sum in terms of the Chern classes of the summands.  But
the formula for the Chern classes of a tensor product of line
bundles is more complicated than usual one. Quillen
showed~\cite{Qui:FGL,Ad:SHGH} that it is as complicated as it can
be.  If $E$ is a complex oriented cohomology theory, then there is
a unique power series
\[
F[s,t]\in \pi_{\ast}E\LL s,t\RR
\]
with the property that for two complex line bundles $L_{1}$ and
$L_{2}$ one has
\[
c_{1}(L_{1}\otimes L_{2})=F[c_{1}(L_{1}),c_{1}(L_{2})].
\]
The power series $F[s,t]$ is a formal group law over $\pi_{\ast}E$.
Quillen showed that when $E=MU$, the resulting formal group law is
universal in the sense that if $F$ is any formal group law over a ring
$R$, then there is a unique ring homomorphism $MU_{\ast}\to R$
classifying $F$.  In this way the complex cobordism spectrum becomes a
topological model for the moduli space of formal group laws.

\subsection{The Adams spectral sequence}

\vskip-5mm \hspace{5mm}

There are exceptions, but for the most part what computations can
be made of the stable and unstable homotopy groups involve
approximating a space by the spaces of a spectrum $E$ whose
homotopy groups are known, or at least qualitatively understood. A
mechanism for doing this was discovered by Adams~\cite{Ad:SS} in
the case $E=H\Z/p$, and later for a general cohomology theory by
both Adams~\cite{adams69:_lectur,adams69:_stabl} and
Novikov~\cite{novikov67:_rings_adams,novikov67:_method}.  The
device is known as the {\em $E$-Adams spectral sequence for $X$},
or, in the case $E=MU$, the {\em Adams-Novikov spectral sequence
for $X$} (or, in case $X=S^{0}$, just the Adams-Novikov spectral
sequence).

The Adams-Novikov spectral sequence has led to many deep insights in
algebraic topology (see, for example,~\cite{Rav:MU,ravenel92:_nilpot}
and the references therein).  It is usually displayed in the first
quadrant, with the groups contributing to $\pi_{k}S^{0}$ all having
$x$-coordinate $k$.  The $y$-coordinate is the $MU$-Adams filtration,
and can be described as follows: a stable map $f:S^{k}\to S^{0}$ has
filtration $\ge s$ if there exists a factorization
\[
S^{k}=X_{0}\to X_{1}\to\dots \to X_{s-1}\to X_{s}=S^{0}
\]
with the property that each of the maps $MU_{\ast}X_{n}\to
MU_{\ast}X_{n+1}$ is zero.  There is a geometric interpretation of
this filtration: a framed manifold $M$ has filtration $\ge s$ if it
occurs as a codimension $n$ corner in a manifold $N$ with corners,
equipped with suitable almost complex structures on its faces
(see~\cite{laures00}).  The Adams-Novikov spectral starts with the
purely algebraic object
\[
E_{2}^{s,t}=\ext^{s,t}_{MU_{\ast}MU}\left(MU_{\ast},MU_{\ast}
\right).
\]
The quotient of the subgroup of $\pi^{k}S^{0}$ consisting of elements
of Adams-Novikov filtration at least $s$, by the subgroup of those of
filtration at least $(s+1)$ is a sub-quotient of the group $\ext^{s,t}$
with $(t-s)=k$.

\subsection{Asymptotics}

\vskip-5mm \hspace{5mm}

For a number $k$, let $g(k)=s$ be the largest integer $s$ for
which $\pi_{k}S^{0}$ has a non-zero element of Adams-Novikov
filtration $s$.  The graph of $g$ is the {\em $MU$-vanishing
curve}, and the main result of~\cite{DHS} is equivalent to the
formula
\[
\lim_{k\to\infty}\frac{g(k)}{k}=0.
\]
This formula encodes quite a bit of the large scale structure of the
category of spectra (see\cite{DHS,HS,ravenel92:_nilpot}), and it would
be very interesting to have a more accurate asymptotic expression.
This is special to complex cobordism.  In the case of the
original Adams spectral sequence for a finite CW complex $X$ (based on
ordinary homology with coefficients in $\Z/p$), it can be
shown~\cite{hopkins99:_vanis_adams} that
\[
\lim_{k\to\infty}\frac{g_{X}^{H\Z/p}(k)}{k}=\frac{1}{2(p^{m}-1)},
\]
for some $m$.  This integer $m$ is an invariant of $X$ know as the
``type'' of $X$.  It coincides with the largest value $m$ for which
the Morava $K$-group $K(m)_{\ast}X$ is non-zero.  For more on the role
of this invariant in algebraic topology, see~\cite{HS,Hop:GMHT}.

Now the $E_{2}$-term of the Adams-Novikov spectral sequence is far
from being zero above the curve $g(n)$, and a good deal of what
happens in spectral sequence has to do with getting rid of what is
up there.  A few years ago, Haynes Miller, and I constructed a
series of spectra designed to classify and capture the way this
happens.  We were motivated by connective $KO$-theory, whose
Adams-Novikov spectral sequence more or less coincides with the
Adams-Novikov spectral for the sphere above a line of slope $1/2$,
and is very easy to understand below that line (and in fact
connective $KO$ can be used to capture everything above a line of
slope $1/5$~\cite{Mah:ImJ,LM,Mah:bo}).  By analogy we called these
cohomology theories $EO_{n}$.  These spectra were used to solve
several problems about the homotopy groups of spheres.

The theory we now call $\tmf$ was originally constructed to
isolate the ``slope $1/6^{\text{th}}$-sector'' of the Adams
Novikov spectral sequence, and in~\cite{hokins95:_topol_witten},
for the reasons mentioned above, it was called $eo_{2}$.  In the
next section the spectrum $\tmf$ will be constructed as a
topological model for the moduli space (stack) of generalized
elliptic curves.

\section{tmf} \label{sec:tmf}

\vskip-5mm \hspace{5mm}

\subsection{The algebraic theory of modular forms}

\vskip-5mm \hspace{5mm}

Let $C$ be the projective plane curve given by the Weierstrass
equation
\begin{equation}\label{eq:1}
y^{2}+a_{1}xy+a_{3}y=x^{3}+a_{2}x^{2}+a_{4}x + a_{6}
\end{equation}
over the ring
\[
A=\Z[a_{1},a_{2},a_{3},a_{4},a_{6}].
\]
Let
\[
A_{\ast}=\bigoplus_{n\in\Z}A_{2n},
\]
be the graded ring with
\[
A_{2n}=H^{0}(C;(\Omega^{1}_{C/A})^{\otimes n}).
\]
If $u\in A_{2}$ is the differential
\[
u=\frac{dx}{2y+a_{1}x+a_{3}},
\]
then
\[
A_{\ast}\approx A[u^{\pm1}].
\]

The $A$-module $A_{2}$ is free over $A$ of rank $1$, and is the module
of sections of the line bundle
\[
\omega := H^{0}\left(\mathcal O_{C}(-e)/\mathcal
O_{C}(-2e) \right) \approx p_{\ast}\Omega^{1}C.
\]
In this expression $p:C\to\spec A$ is the structure map, and $e:\spec
A\to C$ is the point at $\infty$.

Let $G$ be the algebraic group of projective transformations
\begin{align*}
x &\mapsto \lambda^{2}x+r, \\
y &\mapsto \lambda^{3}y + s x + t.
\end{align*}
Such a transformation carries $C$ to the curve $C'$ defined by
an equation
\[
y^{2}+a'_{1}xy+a'_{3}y=x^{3}+a'_{2}x^{2}+a'_{4}x + a'_{6},
\]
for some $a'_{i}$.  This defines an action of $G$ on $A_{\ast}$.  The
ring of invariants
\[
H^{0}(G;A_{\ast})
\]
is the ring of {\em modular forms over $\Z$}.

The structure of $H^{0}(G;A_{\ast})$ was worked out by Tate (see
Deligne~\cite{Deligne:Formulaire}).  After inverting $6$ and
completing the square and the cube, equation~\eqref{eq:1} can be
put in the form
\[
\tilde y^{2} = \tilde x^{3}+\tilde c_{4}\,\tilde x + \tilde c_{6}\qquad \tilde c_{4}, \tilde
c_{6}\in A[\tfrac16]  ,
\]
with
\[
\tilde x=x+\frac{a_{1}^{2}+4 a_{2}}{12}\qquad
\tilde y = y + \frac{a_{1}x+a_{3}}{2}.
\]
The elements
\begin{align*}
c_{4} &= 48\,u^{4}\tilde c_{4}, \\
c_{6} &= 864\,u^{6}\tilde c_{6},
\end{align*}
lie in $A_{\ast}$, and
\[
H^{0}(G;A_{\ast}) =
\Z[c_{4},c_{6},\Delta]/(c_{4}^{3}-c_{6}^{2}=1728\Delta).
\]
We'll write
\[
M_{n}=H^{0}(G;A_{2n})
\]
for the homogeneous part of degree $2n$.  It is the group of {\em
modular forms of weight $n$ over $\Z$}.

\subsection{The topological theory of modular forms}

\vskip-5mm \hspace{5mm}

In~\cite{hopkins95:_topol_witten,Hop:ell1} it is shown that this
algebraic theory refines from rings to {\em ring spectra}, leading
to a topological model for the theory of elliptic curves and
modular forms.  Here is a rough idea of how it goes.

The set of regular points of $C$ has a unique group structure in which
the point at $\infty$ is the identity element, and in which collinear
points sum to zero.  Expanding the group law in terms of the
coordinate $t=x/y$ gives a formal group law
\[
C^{f}[s,t]\in A\LL s,t\RR
\]
over $A$, which, by Quillen's theorem (see
\S\ref{sec:algebr-struct-spectr}) is classified by a graded ring
homomorphism
\[
MU_{\ast}\to A_{\ast}.
\]
The functor
\[
X\mapsto MU_{\ast}(X)\otimes_{MU_{\ast}}A_{\ast}
\]
is not quite a cohomology theory, but it becomes one after inverting $c_{4}$
or $\Delta$.

Based on this, a spectrum $E_{A}$ can be constructed with
\[
\pi_{\ast}E_{A}=A_{\ast},
\]
and representing a complex oriented cohomology theory in which the
formula for the first Chern class of a tensor product of complex line
bundles is given by
\[
c_{1}(L_{1}\otimes L_{2}) = u^{-1}\, C^{f}\left(u\, c_{1}(L_{1}),u\,
c_{1}(L_{2}) \right).
\]
A spectrum $E_{G}$ can be constructed out of the affine coordinate
ring of $G$ in a similar fashion, as can an ``action'' of $E_{G}$ on
$E_{A}$.  The spectrum $\tmf$ is defined to be the $(-1)$-connected
cover of the homotopy fixed point spectrum of this group action.

To actually carry this out requires quite a bit of work.  The
difficulty is that the theory just described only defines an action of
$E_{G}$ on $E_{A}$ up to homotopy, and this isn't rigid enough to form
the homotopy fixed point spectrum.  In the end it can be done, and
there turns out to be an unique way to do it.

\subsection{The ring of topological modular forms}

\vskip-5mm \hspace{5mm}

The spectrum $\tmf$ is a homotopy theoretic refinement of the ring
$H^{0}(G;A_{\ast})$, there is a spectral sequence
\[
H^{s}\left(G;A_{t} \right) \Rightarrow \pi_{t-s}\tmf.
\]
The ring $\pi_{\ast}\tmf$ is the ring of {\em topological
modular forms}, and the group $\pi_{2n}\tmf$ the group of {\em
topological modular forms of weight $n$}.
The edge homomorphism of this spectral sequence is a homomorphism
\[
\pi_{2n}\tmf\to M_{n}.
\]
This map isn't quite surjective, and there is the following
result of myself and Mark Mahowald

\begin{prop}\label{thm:4}
The image of the map $\pi_{2\ast}\tmf\to M_{\ast}$ has a basis given by the
monomials
\[
a_{i,j,k}\,c_{4}^{i}c_{6}^{j}\Delta^{k}\qquad i,k\ge 0, j=0,1
\]
where
\[
a_{i,j,k}=
\begin{cases}
1 &\qquad i>0,j=0 \\
2 &\qquad j=1 \\
24/\gcd(24,k) &\qquad i,j=0.
\end{cases}
\]
\end{prop}

In the table below I have listed the first few homotopy groups of $\tmf$
\bigskip
\begin{center}
\begin{tabular}{|c|c|c|c|c|c|c|c|c|c|}\hline
$k$ & $0$ & $1$ & $2$ & $3$ & $4$ & $5$ & $6$ & $7$ & $8$   \\ \hline
$\pi_{k}\tmf$ & $\Z$ & $\Z/2$ & $\Z/2$  & $\Z/24$ & $0$ & $0$ &
$\Z/2$ & 0 & $\Z \oplus \Z/2$ \\ \hline
\end{tabular}

\vskip .2in

\begin{tabular}{|c|c|c|c|c|c|c|c|c|c|}\hline  $9$ & $10$ & $11$ &
$12$ & $13$ & $14$ & $15$   \\ \hline
$(\Z/2)^{2}$ & $\Z/6$ & $0$ & $\Z$ & $\Z/3$ & $\Z/2$ &
$\Z/2$  \\ \hline
\end{tabular}
\end{center}
\medskip
The homotopy homomorphism induced by the
unit $S^{0}\to\tmf$ of the ring $\tmf$ is the
``$\tmf$-degree,'' a ring homomorphism
\[
\pi_{\ast}S^{0}\to\pi_{\ast}\tmf.
\]
The $\tmf$-degree is an isomorphism in dimensions $\le 6$, and it take
non-zero values on each of the classes represented by the Lie groups
$SU(3)$, $U(3)$, $Sp(2)$, $Sp(1)\times Sp(2)$, $G_{2}$, $U(1)\times
G_{2}$, regarded as framed manifolds via their left invariant
framings.  Thus, combined with the Hopf-invariant and the invariants
coming from $KO$-theory, the $\tmf$-degree accounts for all of
$\pi_{\ast}S^{0}$ for $\ast\le 15.$ In fact the Hopf-invariant and the
invariants coming from $KO$ can also be described in terms of $\tmf$
and nearly all of $\pi_{\ast}S^{0}$ for $\ast<60$ can be accounted
for.

\section{{\boldmath $\theta$}-series}\label{sec:vartheta-series}

\vskip-5mm \hspace{5mm}

\subsection{Cohomology rings as rings of functions}
\label{sec:cohomology-rings-as}

\vskip-5mm \hspace{5mm}

Consider the computation
\[
H^{\ast}(\cpinfty;\Z)=\Z[x].
\]
On one hand this tells us something about the cell structure of
complex projective space;   the cohomology class $x^{n}$ is ``dual''
to the cell in dimension $2n$.  On the other hand, a polynomial is a
function on the affine line, and the elements of $H^{\ast}(\cpinfty)$
tell us something about the affine line.  Combining these, the
prospect presents itself, of using the cell structure of one space to
get information about the function theory of another.

We will apply this not to ordinary cohomology, but to the
cohomology theory $E_{A}$.  Before doing so, more of the function
theoretic aspects of $E_{A}$ need to be spelled out.  By
construction, the ring $E_{A}^{0}(\cpinfty)$ is the ring of
functions on the formal completion of $C$ at the point $e$ at
$\infty$.  The ring $E_{A}^{0}(\hpinfty)$ is the ring of functions
which are invariant under the involution
\begin{equation}\label{eq:7}
\begin{aligned}
\tau(x) &= x \\
\tau(y) &= -y-a_{1}x-a_{3}
\end{aligned}
\end{equation}
given by the ``inverse'' in the group law.  The map
\[
E_{A}^{0}(\cpinfty)\to E_{A}^{0}(\text{pt})
\]
corresponds to evaluation at $e$, and the reduced cohomology group
$\tilde E_{A}^{0}(\cpinfty)$ to the ideal of formal function vanishing
at $e$.  Note that this is consistent with the definition of $A_{2}$
as sections of the line bundle $\omega$:
\[
A_{2}=\pi_{2}E_{A}=\tilde E_{A}^{0}(S^{2})=I/I^{2}=H^{0}\left(\mathcal
O_{C}(-e)/\mathcal O_{C}(-2e) \right).
\]
Now $\tilde E_{A}^{0}(\cpinfty)$ is the cohomology group of the Thom
complex $E_{A}^{0}({\cpinfty}^{L})$.  More generally, there is an
additive correspondence
\[
\left\{\text{virtual representations of
$U(1)$}\right\}\leftrightarrow \left\{\text{divisors on $C$}
\right\},
\]
under which a virtual representation $V$ corresponds to a divisor
$D$ for which
\[
E^{0}\left((\cpinfty)^{V} \right)=H^{0}\left(C^{f}; {\mathcal
0}(D)\otimes\left(\Omega^{1}\right)^{\otimes\mu(V)} \right),
\]
where $\mu(V)$ is the multiplicity of the trivial representation
in $V$.  There is a similar correspondence between even functions
with divisors of the form $D+\tau^{\ast}D$, virtual
representations $V$ of $SU(2)$ and $E_{A}^{0}\left((\hpinfty)^{V}
\right)$ (with $\tau$ the involution~\eqref{eq:7}).

\subsection{The Hopf fibration and the Weierstrass {\boldmath $\pfn$}-function}
\label{sec:hopf-fibr-weierstr}

\vskip-5mm \hspace{5mm}

Consider the function $x$ in the Weierstrass
equation~\eqref{eq:1}. This function has a double pole at $e$.  We
now ask if there is a ``best'' $x$ to choose, ie, a function with
a double pole at $e$ which is invariant under the action of $r$,
$s$ and $t$.  Such a function will be an eigenvector for $\lambda$
with eigenvalue $\lambda^{2}$. It is more convenient to search for
an quantity which is invariant under $\lambda$ as well, so instead
we search for a quadratic differential on $C$, i.e. a section
\[
x'\in H^{0}\left(C;\Omega^{1}(e)^{\otimes 2} \right)
\]
which is invariant under $r,s,t$ and $\lambda$.  Now the space
$H^{0}(C;\Omega^{1}(e)^{2})$ has dimension $2$, and sits in a
short exact sequence of vector spaces over $\spec A$
\begin{equation}\label{eq:3}
0\to\omega^{2}\to H^{0}(C;\Omega^{1}(e)^{2})\to \mathcal O_{A}\to 0,
\end{equation}
where $\omega$ is the line bundle of invariant differentials on $C$,
and the second map is the ``residue at $e$''.  This sequence is
$G$-equivariant, and defines an element of
\[
\ext^{1}(\mathcal O, \omega^{2})=H^{1}(G;A_{4}).
\]
The obstruction to the existence of an $x'$ with residue $1$ is the
Yoneda class $\nu$ of this extension.  Completing the square and cube
in~\eqref{eq:1}, gives the $G$-invariant expression
\[
\left(12 x + a_{1}^{2}+ 4 a_{2} \right) u^{2},
\]
so that $12\nu=0$.  In fact the group $H^{1}(G;A_{4})$ is cyclic of
order $12$ with $\nu$ as a generator.  The group $\pi_{3}\tmf=\Z/24$ is
assembled from $H^{1}(G;A_{4})$ and $H^{3}(G;A_{6})=Z/2$, and sits in
an exact sequence
\[
0\to H^{3}(G;A_{6})\to  \pi_{3}\tmf\to H^{1}(G;A_{4})\to 0.
\]

This $12$ can also be seen transcendentally.  Over the complex
numbers, a choice of $x$ is given by the Weierstrass $\pfn$-function:
\[
\pfn(z,\tau)=
\frac{1}{z^{2}}+\sum_{0\ne(m,n)\in\Z^{2}}
\frac{1}{(z-m\tau-n)^{2}}-\frac{1}{(m\tau+n)^{2}}.
\]
The Fourier expansion of $\pfn(z,\tau)\,dz^{2}$ is
(with $q=e^{2\pi\,i\,\tau}$ and $u=e^{2\pi\,i\,z}$)
\[
\pfn(z,\tau)\,dz^{2}=\left(\sum_{n\in\Z}
\frac{q^{n}u}{(1-q^{n}u)^{2}}+\frac{1}{12}
-2\sum_{n\ge 1}\frac{q^{n}}{(1-q^{n})^{2}} \right)
\left(\frac{du}{u} \right)^{2}.
\]
Note that all of the Fourier coefficients of $\pfn$ are integers, except
for the constant term, which is $1/12$.  This is the same $12$.

Under the correspondence between divisors and Thom complexes, the
differential $x'$ corresponds to a $G$-invariant element
\[
x'_{\text{top}}\in E_{A}^{0}\left(\hp^{(2-V)} \right),
\]
with $V$ the defining representation of $SU(2)$.  Now the spectrum
$\hp^{2-V}$ has a (stable) cell decomposition
\[
\hp^{2(1-L)}=S^{0}\cup_{\nu} e^{4}\cup\dots
\]
with one cell in every real dimension $4k$.  The $4$-cell is attached
to the $0$-cell by the stable Hopf map $\nu:S^{7}\to S^{4}$, which
generates $\pi_{3}(S^{0})=\Z/24$.  The restriction of the quadratic
differential $x'_{\text{top}}$ to the zero cell is given by the
residue at $e$, and the obstruction to the existence of an
$G$-invariant $x'_{\text{top}}$ with residue $k$ is the image of $k$
under the connecting homomorphism
\begin{equation}
\label{eq:5}
H^{0}(G;A)\to H^{1}\left(G;E_{A}^{0}\left(
\hp^{2-V}/S^{0}\right)\right).
\end{equation}

To evaluate~\eqref{eq:5}, note that the map
\[
H^{1}\left(G;E_{A}^{0}\left(
\hp^{2-V}/S^{0}\right)\right)\to
H^{1}\left(G;E_{A}^{0}(S^{4}) \right)=
H^{1}\left(G;A_{4}) \right)=\Z/12
\]
is projection onto a summand, and the image of~\eqref{eq:5} is
contained in this summand.    Thus the obstruction to the existence of
the quadratic differential $x_{\text{top}}$ is the same $k\in\Z/12$.
In this way, the theory of topological modular forms relates the Hopf
map $\nu$ to the constant term in the Fourier expansion of the
Weierstrass $\pfn$-function, and to the existence of a certain
quadratic differential on the universal elliptic curve.

\subsection{Lattices and their {\boldmath $\theta$}-series}

\vskip-5mm \hspace{5mm}

There is a slightly more sophisticated application of these ideas
to the theory of even unimodular lattices.  Suppose that $L$ is a
positive definite, even unimodular lattice of dimension $2d$.  The
{\em theta function} of $L$, $\theta_L$ is the generating function
\begin{align*}
\theta_L(q) & =\sum_{\ell\in L} q^{\frac12\langle \ell,\ell\rangle}\\
&= \sum_{n\ge 0} L_n q^n, \\
L_n &= \#\{\ell\mid \langle l,l\rangle = 2n \}.
\end{align*}
It follows from the Poisson summation formula that $\theta_L(q)$ is
the $q$-expansion of a modular form over $\Z$ of weight $d$, and so
lies in the ring
\begin{align*}
\Z[c_4,c_6,\Delta]/&(c_4^3-c_6^2-1728\Delta)  \subset Z\LL q\RR ,\\
c_4 &= 1 + 240 \sum_{n>0} \sigma_3(n) q^n ,\\
c_6 &= 1 - 504\sum_{n>0} \sigma_5(n) q^n ,\\
\Delta &= q\prod_{n=1}^{\infty} (1-q^n)^{24}.
\end{align*}
Since the group of modular forms of a given weight is finitely
generated, the first few $L_n$ determine the rest.  This leads to
many restrictions on the distributions of lengths of vectors in a
positive definite, even unimodular lattice.

The $\theta$-series of $L$ is the value at $z=0$ of the $\theta$-function
\[
\theta(z,\tau)=\sum_{\ell\in L} e^{\pi i \langle\ell,\ell
\rangle\tau+2\pi i\langle\ell,z \rangle}, \qquad z\in \C\otimes
L,\quad q=e^{2\pi i\, \tau},
\]
which, under the correspondence between divisors and representations
has the following topological interpretation.  Let $V$ be any
$d$-dimensional (complex) virtual representation of $U(1)\otimes L$
with the property that $c_{1}(V)=0$ and $c_{2}(V)$ corresponds to the
quadratic form, under the isomorphism
\[
H^{2}\left(B(U(1)\otimes L);\Z\right)=\sym^{2}L^{\ast}.
\]
For example, if $A=(a_{ij})$ is the matrix of the quadratic form with
respect to some basis, then $V$ could be taken to be
\[
V= \C^{d} + \frac12\sum_{i,j}a_{ij}(1-L_{i})(1-L_{j}),
\]
where $L_{i}$ is the character of $U(1)\otimes L$ dual to the
$i^{\text{th}}$ basis element.  The series $\theta(z,\tau)$
corresponds to a $G$-invariant element
\[
\theta^{\text{top}}\in E_{A}^{0}(\left(BU(1)\otimes L \right)^{V}).
\]
The restriction of $\theta^{\text{top}}$ to $\left\{\text{pt}
\right\}^{V}=S^{2n}$ is an element of $H^{0}(G;A_{2d})$, i.e. an
algebraic modular form of weight $d$.  This modular form is
$\theta_{L}$.

Now the Thom spectrum $\left(BU(1)\otimes L \right)^{V}$ has a
stable cell decomposition
\[
S^{2d}\cup \bigvee^{2d} e^{2d+2}\cup \bigvee^{d(2d+1)}
e^{2d+4}\cup\cdots.
\]
Since $c_{1}(V)=0$, the cells of dimension $2d+2$ are not attached
to the $(2d)$-cell.  The assumptions on the quadratic form, and on
$c_{2}(V)$ imply that one of the $(2d+4)$-cells is attached to the
$(2d)$-cell by (a suspension of) the stable Hopf map $\nu$.  The
presence of this attaching map implies the following mod $24$
congruence on $\theta_{L}$.
\begin{thm}\label{thm:3}
Suppose $L$ is a positive definite, even unimodular lattice of
dimension $24k$.
Write
\[
\theta_L(q) = c_4^{3k} + x_1 c_4^{3(k-1)}\Delta + \dots + x_k\Delta^k.
\]
Then
\begin{equation}\label{eq-borcherds-delta}
x_k\equiv 0 \mod 24.
\end{equation}
\end{thm}

The above result was originally proved by
Borcherds~\cite{borcherds95:_autom_o_r} as part of his investigation
into infinite product expansions for automorphic forms on certain
indefinite orthogonal groups.  The above topological proof can be
translated into the language of complex function theory.  The details
are in the next section.

The congruence of Theorem~\ref{thm:3} together with
Proposition~\ref{thm:4} give the following
\begin{prop}
\label{thm:2}
Suppose $L$ is a positive definite, even unimodular lattice of
dimension $2d$.  There is an element $\theta_{L}^{\text{top}}\in
\tmf^{0}(S^{2d})$ whose image in $M_{d}$ is $\theta_{L}$.
\end{prop}

It can also be shown that the $G$-invariant
\[
\theta^{\text{top}}\in H^{0}(G;E_{A}^{0}(B(U(1)\otimes L)^{V}))
\]
is truly topological in the sense that it is the representative in
the $E_{2}$-term of the spectral sequence
\[
H^{0}(G;E_{A}^{0}(B(U(1)\otimes L)^{V})) \Rightarrow
\tmf^{0}(B(U(1)\otimes L)^{V}),
\]
of an element in $\tmf^{0}(B(U(1)\otimes L)^{V})$.  I don't know of a
direct construction of these truly topological theta series.

\subsection{An analytic proof of Theorem~\ref{thm:3}}

\vskip-5mm \hspace{5mm}

The analytic interpretation of the proof of Theorem~\ref{thm:3}
establishes the result in the form
\begin{equation}\label{eq-borcherds-residue}
\res_{\Delta=0}\frac{\theta_L(q)}{\Delta^{k}}\, \frac{d\Delta}{\Delta} \equiv
\res_{q=0}\frac{\theta_L(q)}{\Delta^{k}}\, \frac{dq}{q} \equiv 0 \mod 24.
\end{equation}
The equivalence of~\eqref{eq-borcherds-delta}
with~\eqref{eq-borcherds-residue} follows easily from the facts
\begin{align*}
c_4 &\equiv 0 \mod 24 ,\\
\frac{d\Delta}{\Delta} &\equiv \frac{dq}{q} \mod 24.
\end{align*}

Set $u=e^{2\pi i z}$, and for a vector $\mu\in L$ let
\[
\phi_\mu(z,\tau) = \frac{\sum_{\ell\in L}q^{\frac12\langle
\ell,\ell\rangle}u^{\langle \mu,\ell\rangle}}
{\sigma(q,u)^{\langle \mu,\mu\rangle}},
\]
where
\[
\sigma(q,u) = u^{\frac12}(1-u^{-1})
\prod_{n=1}^{\infty}\frac{(1-q^n u)(1-q^nu^{-1})}
{(1-q^n)^2}
\]
is the Weierstrass $\sigma$-function.  It is immediate
from the definition that
\[
\phi_{\mu}(-z,\tau) = \phi_{\mu}(z,\tau),
\]
and it follows from the modular transformation formula for $\theta$ that for
\[
m,n\in\Z, \qquad
\begin{pmatrix}
a & b\\
c& d
\end{pmatrix}\in \SL_2(Z),
\]
the function
$\phi_{\mu}(z,\tau)$ satisfies
\begin{align}
\phi_{\mu}(z + m\tau+n,\tau) &=\phi_{\mu}(z,\tau),\\
\label{eq-modular-phi}
\phi_{\mu}\left(\frac z{c\tau +d},\frac{a\tau +
b}{c\tau+d}\right) & = (c\tau+d)^{d/2}\,\phi_{\mu}(z,\tau).
\end{align}
These identities are equivalent to saying that if we write,
with $x=2\pi i z$,
\[
\phi_{\mu}(z,\tau) = x^{-\langle \mu,\mu
\rangle}\left(\phi_{\mu}^{(0)} + \phi_{\mu}^{(2)} x^2 + \dots\right),
\]
then $\phi_{\mu}^{(2k)}$ is a modular form of weight $d/2 + 2k$.

It is the term $\phi_{\mu}^{(2)}$ that stores the information about
the $(2d+4)$-cells of the spectrum $B(U(1)\otimes L)^{V}$.  A little
computation shows that
\begin{equation}\label{eq-phi}
\phi_{\mu}^{(2)} = \theta_{\mu}-\frac{\langle \mu,\mu\rangle}{24}\theta_L,
\end{equation}
where
\[
\theta_{\mu} = \sum_{\ell\in L} q^{\frac12\langle
\ell,\ell\rangle}\frac{\ellmu^2}{2} - \mumu \sum_{n\ge 1}
\frac{q^n}{(1-q^n)^2}.
\]
Because of~\eqref{eq-modular-phi}, the expression
\[
\frac{\phi_{\mu}^{(2)}}{\Delta^k}\, \frac{dq}{q}
\]
is a meromorphic differential on the projective $j$-line.  It's
only pole is at $q=0$.  It follows that\footnote{I learned this
trick from Borcherds, who told me he learned it from a referee.}
\[
\res_{q=0}\frac{\phi_{\mu}^{(2)}}{\Delta^k}\,\frac{dq}{q} = 0,
\]
and hence
\begin{equation}\label{eq-main-residue}
\res_{q=0} \frac{\theta_{\mu}}{\Delta^k}\frac{dq}{q}
=
\frac{\langle \mu,\mu\rangle}{24}\res_{q=0}\frac{\theta_L}{\Delta^k}\frac{dq}{q}.
\end{equation}

To deduce~\eqref{eq-borcherds-residue}, define
\[
p:L\otimes \Z/2 \to \Z/2
\]
by letting $p(\mu)$ be the value of~\eqref{eq-main-residue} reduced
modulo $2$.  The left hand side of~\eqref{eq-main-residue} shows that
$p$ does in fact takes its values in $\Z/2$.  Making use of the
symmetry $\ell \mapsto -\ell$, it also shows that $p$ is linear.  The
right hand side shows that $p$ is quadratic with underlying bilinear
form
\[
p(\mu_1+\mu_2)
-p(\mu_1)-p(\mu_2)=
\langle\mu_1,\mu_2\rangle
\frac1{12}
\res_{q=0}\frac{\theta_L}{\Delta^k}\frac{dq}{q}.
\]
It follows that
\[
\langle\mu_1,\mu_2\rangle
\frac1{12}
\res_{q=0}\frac{\theta_L}{\Delta^k}\frac{dq}{q}\equiv 0 \mod 2.
\]
Since $L$ is unimodular, there are vectors $\mu_1$ and $\mu_2$
with
\[
\langle\mu_1,\mu_2\rangle=1
\]
and so
\[
\frac1{12}
\res_{q=0}\frac{\theta_L}{\Delta^k}\frac{dq}{q}\equiv 0 \mod 2.
\]
This is~\eqref{eq-borcherds-residue}.

\section{Topological modular forms as cobordism invariants}
\label{sec:topol-modul-forms}

\vskip-5mm \hspace{5mm}

\subsection{The Atiyah-Bott-Shapiro map}
\label{sec:atiyah-bott-shapiro}

\vskip-5mm \hspace{5mm}

As described in \S 2, the part of the homotopy groups of spheres
that is best understood geometrically is the part that is captured
by $KO$-theory.  The geometric interpretations rely on the fact
that all of the corresponding framed cobordism invariants can be
expressed in terms of invariants of $\spin$-cobordism. From the
point of view of topology what makes this possible is the
factorization
\[
S^{0}\to\mspin\to KO
\]
of the unit in $KO$-theory, through the Atiyah-Bott-Shapiro map
$\mspin\to KO$.  The Atiyah-Bott-Shapiro map is constructed using the
representations of the spinor groups, and relies on knowing that for
a space $X$, elements of $KO^{0}(X)$ are represented by vector bundles
over $X$.  It gives a $KO$-theory Thom isomorphism for $\spin$-vector
bundles, and is topological expression for the index of the Dirac
operator.  Until recently it was not known how to produce this map by
purely homotopy theoretic means.

The framed cobordism invariants coming from $\tmf$ cannot be expressed
directly in terms of $\spin$-cobordism.  One way of seeing this is to
note that the group $\pi_{3}\mspin$ is zero.  The fact that the map
\[
\pi_{3}S^{0}\to\pi_{3}\tmf
\]
is an isomorphism prevents a factorization of the form
\[
S^{0}\to\mspin\to\tmf.
\]
In some sense this is all that goes wrong.  Let $\boeight$ be the
$7$-connected cover of $\bspin$, and $\moeight$ the corresponding Thom
spectrum.  We will see below that there is a factorization of the unit
\[
S^{0}\to\moeight\to\tmf.
\]

There is not yet a geometric interpretation of $\tmf^{0}(X)$, so the
construction of a map $\moeight\to\tmf$ must be made using homotopy
theoretic methods.  The key to doing this is to exploit the
$\einfty$-ring structures on the spectra involved.

\subsection{{\boldmath $\einfty$}-maps}

\vskip-5mm \hspace{5mm}

The spectra $\moeight$, $\mspin$, $\tmf$ and $KO$ are all
$\einfty$-ring spectra, and from the point of view of homotopy
theory it turns out easier to construct $\einfty$ maps
\[
\mspin\to KO\quad\text{and}\quad\moeight\to \tmf
\]
than merely to produce maps of spectra.  In fact the homotopy type
of the spaces of $\einfty$-maps
\begin{equation*}
\einfty\left(\moeight,\tmf \right) \quad\text{and}\quad
\einfty\left(\mspin,KO \right)
\end{equation*}
can be fairly easily identified using homotopy theoretic methods.
In this section I will describe their sets of path components.

A map $\phi:\mspin\to KO$ determines a Hirzebruch genus with an even
characteristic series $K_{\phi}(z)\in \Q\LL x\RR$.  Define the {\em
characteristic sequence} of $\phi$ to be the sequence of rational
numbers
\[
(b_{2},b_{4},\cdots)
\]
given by
\[
\log\left(K_{\phi}(z) \right)= -2\sum_{n>0}  b_{n}\, \frac{x^{n}}{n!}.
\]
We use this sequence to form the {\em characteristic map}
\[
\pi_{0}\einfty\left(\mspin,KO \right)\to
\left\{(b_{2},b_{4},\cdots)\mid b_{2i}\in\Q \right\}
\]
from the set of homotopy classes of $\einfty$ maps $\mspin\to KO$ to
the set of sequences of rational numbers.

Let $B_{n}$ denote the $n^{\text{th}}$ Bernoulli number defined by
\[
\frac{x}{e^{x}-1} = \sum_{n\ge 0} B_{n}\frac{x^{n}}{n!}.
\]
The following result is due to Matthew Ando, myself, and Charles Rezk:

\begin{thm}
\label{thm:8}
The characteristic map
\[
\pi_{0}\einfty\left(\mspin,KO \right)\to
\left\{(b_{2},b_{4},\cdots)\mid b_{2i}\in\Q \right\}
\]
gives an isomorphism of $\pi_{0}\einfty\left(\mspin,KO \right)$ with
the set of sequences
\[
(b_{2},b_{4},b_{6},\dotsc )
\]
for which
\begin{thmList}
\item $b_{n}\equiv B_{n}/2n\mod\Z$;
\item for each odd prime $p$ and each $p$-adic unit $c$,
\begin{multline*}
m \equiv n \mod p^{k}(p-1) \\
\implies (1-c^{n})(1-p^{n-1})b_{n}\equiv
(1-c^{m})(1-p^{m-1})b_{m}\mod p^{k+1};
\end{multline*}
\item for each $2$-adic unit $c$
\begin{multline*}
m \equiv n \mod 2^{k} \\
\implies (1-c^{n})(1-2^{n-1})b_{n}\equiv
(1-c^{m})(1-2^{m-1})b_{m}\mod 2^{k+2}.
\end{multline*}
\end{thmList}
\end{thm}

\begin{rem}\rm
By the Kummer congruences, the sequence with $b_{2n}=B_{2n}/(4n)$
comes from an $\einfty$-map $\mspin\to KO$.  The associated
characteristic series is
\[
\frac{x}{e^{x/2}-e^{-x/2}}=\frac{x/2}{\sinh(x/2)}
\]
and the underlying map of spectra coincides with the one constructed
by Atiyah-Bott-Shapiro.  Theorem~\ref{thm:8} therefore gives a purely
homotopy theoretic construction of this map.
\end{rem}

We now turn to describing the set of homotopy classes of $\einfty$
maps
\[
\moeight\to\tmf.
\]
Associated to a multiplicative map $\moeight\to\tmf$ is a
characteristic series of the form
\[
K_{\phi}(z)=\sum a_{2n} z^{2n}\qquad a_{2n}\in \Q\LL q\RR
\]
which is well-defined up to multiplication by the exponential of a
quadratic function in $z$.  We associate to such a series, a sequence
\begin{equation}
\label{eq:8}
(g_{4},g_{6}\cdots)\quad g_{2n}\in \Q\LL q\RR
\end{equation}
according to the rule
\[
\log\left(K_{\phi}(z) \right)= -2\sum_{n>0}  g_{n}\, \frac{x^{n}}{n!}.
\]
For $n>2$ the terms $g_{n}$ is independent of the quadratic
exponential factor, and is the $q$-expansion of a modular form of
weight $n$ (see~\cite{hopkins95:_topol_witten, ando01:_ellip_witten}).
This defines the {\em characteristic map}
\[
\pi_{0}\einfty\left(\moeight,\tmf \right)\to
\left\{(g_{2},g_{4},\cdots) \mid g_{2n}\in M_{2n}\otimes\Q\right\}.
\]

Following Serre~\cite{serre73:_formes} let $G_{2k}$ denote the
(un-normalized) Eisenstein series of weight $2k$
\[
G_{2k}=-\frac{B_{2k}}{4k}+\sum_{n>0}\sigma_{2k-1}(n)q^{n},\qquad
\sigma_{2k-1}(n)=\sum_{d\vert n}d^{2k-1},
\]
and for a prime $p$, let $G_{2k}^{\ast}$ be the (un-normalized) $p$-adic
Eisenstein series
\[
G^{\ast}_{2k} = -(1-p^{2k-1})
\frac{B_{2k}}{4k}+\sum_{n>0}\sigma^{\ast}_{2k-1}(n)q^{n}, \qquad
\sigma^{\ast}_{2k-1}(n)=\sum_{\substack{d\vert n\\(d,p)=1}}d^{2k-1}.
\]
We will also need the Atkin ($U$) and Verschiebung ($V$) operators
on $p$-adic modular forms of weight $k$
(See~\cite[\S2.1]{serre73:_formes}).  For a $p$-adic modular form
\[
f=\sum_{n=0}^{\infty}a_{n}q^{n}
\]
of weight $k$, one defines
\[
f\vert U = \sum_{n=0}^{\infty}a_{pn}q^{n}\qquad
f\vert V = \sum_{n=0}^{\infty}a_{n}q^{pn}.
\]
Finally, we set
\[
f^{\ast}=f-p^{k-1}f\vert V.
\]
(This gives two meanings to the symbol $G^{\ast}_{2k}$, which are
easily checked to coincide.)

\begin{prop}
\label{thm:5}
The image of the characteristic map
\[
\pi_{0}\einfty\left(\moeight,\tmf \right)\to
\left\{(g_{2},g_{4},\cdots) \mid g_{2n}\in M_{2n}\otimes\Q\right\}
\]
is the set of sequences $(g_{2n})$ satisfying
\begin{thmList}
\item $g_{2n}\equiv G_{2n}\mod \Z$
\item For each odd prime $p$ and each $p$-adic unit $c$,
\[
m \equiv n \mod p^{k}(p-1)
\implies (1-c^{n})g^{\ast}_{n}\equiv
(1-c^{m})g^{\ast}_{m}\mod p^{k+1},
\]
\item For each $2$-adic unit $c$,
\[
m \equiv n \mod 2^{k}
\implies (1-c^{n})g^{\ast}_{n}\equiv
(1-c^{m})g^{\ast}_{m}\mod 2^{k+2}.
\]
\item For each prime $p$, $g^{\ast}_{m}\vert U=g^{\ast}_{m}$.
\end{thmList}
The characteristic map is a principle $A$-bundle over its image, where
$A$ is a group isomorphic to a countably infinite product of $\Z/2$'s.
\end{prop}

\begin{rem}\rm
The group $A$ occurring in the final assertion of
Proposition~\ref{thm:5} can be described explicitly in terms of
modular forms.  The description is somewhat technical, and has been
omitted.
\end{rem}

By the Kummer congruences, the sequence of Eisenstein series $(G_{4},
G_{6},\cdots)$ satisfy the conditions of Proposition~\ref{thm:5}.  The
corresponding characteristic series can be taken to be
\[
\frac{z/2}{\sinh(z/2)}
\prod_{n\ge1}\frac{(1-q^n)^2}{(1-q^ne^z )(1-q^ne^{-z})},
\]
and the genus is the Witten genus
$\phi_{W}$~\cite{Wit:Loop,Segal:ell,Hirz:MMF}.  This gives the
following corollary, which was the main conjecture
of~\cite{hokins95:_topol_witten}.

\begin{cor}
\label{thm:6}
There is an $\einfty$-map $\moeight\to\tmf$ whose underlying genus
is the Witten genus.
\end{cor}

\begin{rem}\rm
This refined Witten genus is only specified up to action by an element
of the group $A$ occurring in Proposition~\ref{thm:5}.  It looks as if
a more careful analysis could lead to specifying a single
$\einfty$-map $\moeight\to\tmf$, but this has not yet been carried
out.
\end{rem}

\subsection{The image of the cobordism invariant}

\vskip-5mm \hspace{5mm}

The existence of a refined Witten genus has an application to the
theory of even unimodular lattices.  The following two results are
due to myself and Mark Mahowald.

\begin{thm}
Let $\moeight\to\tmf$ be any multiplicative map whose underlying genus
is the Witten genus.  Then the induced map of homotopy groups
$\pi_{\ast}\moeight\to\pi_{\ast}\tmf$ is surjective.
\end{thm}

Combining this with Proposition~\ref{thm:2} then gives

\begin{cor}\label{thm:7}
Let $L$ be a positive definite, even unimodular lattice of
dimension $2d$.  There exists a $7$-connected manifold $M_{L}$ of
dimension $2d$, whose Witten genus is the $\theta$-function of
$L$, i.e.
\[
\phi_{W}(M)=\theta_{L}.
\]
\end{cor}

In case $L$ is the Leech lattice, the existence of $M_{L}$ gives an
affirmative answer to Hirzebruch's ``Prize Question''~\cite{Hirz:MMF}.

\subsection{Spectra of units and {\boldmath $\einfty$}-maps}
\label{sec:unit-spectra-einfty}

\vskip-5mm \hspace{5mm}

The proofs of Theorems~\ref{thm:8} and~\ref{thm:5} come down to
understanding the structure of the group of units in $KO^{0}(X)$
and $\tmf^{0}(X)$.

For a ring spectrum $R$, let $\Gl(R)$ be the classifying space for
the group of units in $R$-cohomology:
\[
[X,\Gl(R)]=R^{0}(X)^{\times}.
\]
If $R=\left\{R_{n},t_{n} \right\}$, then $\Gl(R)$ is part of the
homotopy pullback square
\begin{equation}
\label{eq:11}
\begin{CD}
\Gl(R) @>>> R_{0}\\
@VVV @VVV \\
\pi_{0}R^{\times}@>>> \pi_{0}R.
\end{CD}
\end{equation}
When $R$ is an $\ainfty$-ring spectrum, $\Gl(R)$ has a classifying
space $\BGl(R)$.  When $R$ is $\einfty$, then $\Gl(R)$ is an infinite
loop space, i.e. there is a spectrum $\gl(R)$ with
$\gl(R)_{0}=\Gl(R)$ (and $\gl(R)_{1}=\BGl(R)$).  The space $\BGl(S^{0})$
is the classifying space for unoriented stable spherical fibrations,
and the map
\begin{equation}
\label{eq:10}
\boeight\to \BGl(S^{0})
\end{equation}
whose associated Thom spectrum is $\moeight$, is an infinite loop
map.   More specifically, let $\sboeight$ be the
$7$-connective cover of the spectrum $KO$.  Then
\[
\left(\Sigma^{-1}\sboeight \right)_{1}=\boeight,
\]
and there is a map of spectra
\[
\Sigma^{7}\sboeight\to\gl(S^{0})
\]
for which the induced map $\left(\Sigma^{-1}\sboeight\right)_{1}\to
\left(\gl(S^{0}) \right)_{1}$ becomes~\eqref{eq:10}.

The following result is what makes it easier to construct
$\einfty$-maps than merely maps of spectra.

\begin{prop}\label{thm:1}
The space $\einfty(\mspin,KO)$ is canonically homotopy equivalent
to the space of factorizations
\[
\xymatrix{
\gl(S^{0}) \ar[r]\ar[d] & \gl(S^{0})\cup
C\Sigma^{-1}\bo\langle4\rangle
 \ar@{-->}[dl] \\
\gl(KO),
}
\]
and the space $\einfty(\moeight,\tmf)$ is canonically homotopy equivalent
to the space of factorizations
\[
\xymatrix{
\gl(S^{0}) \ar[r]\ar[d] & \gl(S^{0})\cup C\Sigma^{-1}\sboeight \ar@{-->}[dl] \\
\gl(\tmf).
}
\]
\end{prop}

\subsection{The Atkin operator and the spectrum of units in {\boldmath $\tmf$}}

\vskip-5mm \hspace{5mm}

Proposition~\ref{thm:1} emphasizes the important role played by
the spectrum of units $\gl(\tmf)$ and $\gl(KO)$.  Getting at the
homotopy type of these spectra uses work of
Bousfield~\cite{bousfield87:_uniquen_k} and
Kuhn~\cite{kuhn89:_morav_k}.  Fix a prime $p$ let $K(n)$ denote
the $n^{\text{th}}$ Morava $K$-theory at $p$, and
\[
L_{K(n)}:\text{Spectra}\to \text{Spectra}
\]
the localization functor with respect to $K(n)$
(see~\cite{Rav:Loc,ravenel92:_nilpot}).  Bousfield (in case $n=1$) and
Kuhn (in case $n>1$) construct a functor
\[
\bkn:\text{Spaces}\to\text{Spectra}
\]
and a natural equivalence of $\bkn(E_{0})$ with $L_{K(n)}E$.  The
spectrum $\bkn(X)$ depends only on the connected component of $X$
containing the basepoint.  In the special case of an $\einfty$-ring
spectrum $E$ it gives (because of~\eqref{eq:11}) a canonical
equivalence
\[
L_{K(n)}\gl E\approx L_{K(n)}E,
\]
which, when composed with the localization map $\gl E\to L_{K(n)}\gl
E$, leads to a ``logarithm''
\[
\log_{K(n)}^{E}:\gl E\to L_{K(n)}E.
\]

Bousfield showed that for $KO$, the logarithm
\[
\log_{K(1)}^{KO}:\gl KO\to L_{K(1)}KO=KO_{p}
\]
becomes an equivalence after completing at $p$ and passing to
$2$-connected covers.  Charles Rezk has recently shown that for a
space $X$, the map
\[
\log_{K(1)}^{KO}:KO^{0}(X)^{\times}\to KO^{0}_{p}(X)
\]
is given by the formula
\[
\frac1p\log\left(\frac{\psi_{p}(x)}{x^{p}}\right).
\]
This equivalence between the multiplicative and additive groups of
$K$-theory was originally observed by Sullivan, and proved by
Atiyah-Segal~\cite{atiyah71:_expon}.

In the case of $\tmf$, Paul Goerss, Charles Rezk and I have shown that
the Bousfield-Kuhn logarithms leads to a commutative diagram
\begin{equation}
\label{eq:9}
\begin{CD}
\gl \tmf @>\log^{\tmf}_{p}>> \tmf_{p} \\
@VVV @VVV \\
L_{K(1)}\tmf @>> 1-U > L_{K(1)}\tmf
\end{CD}
\end{equation}
which becomes homotopy Cartesian after completing at $p$ and passing
to $3$-connected covers.  The spectrum $\tmf_{p}$ is the $p$-adic
completion of $\tmf$, and the map $\log_{p}^{\tmf}$ has a description
in terms of modular forms, but to describe it would take us outside
the scope of this paper.  The spectrum $L_{K(1)}\tmf$ is the
topological analogue of the theory of $p$-adic modular forms of
Serre~\cite{serre73:_formes} and Katz~\cite{katz73}, and the map $U$
is the topological Atkin operator.  One noteworthy feature of this
square is that it locates the Atkin operator in the theory of all
modular forms (and not just $p$-adic modular forms).  Another is that
it connects the failure of $\log_{p}^{\tmf}$ to be an isomorphism with
the spectrum of $U$.  For example, let $F$ denote the fiber of
\[
\gl\tmf_{p} \xrightarrow{\log_{p}^{\tmf}}{} \tmf_{p}.
\]
By the square~\eqref{eq:9}
\[
\pi_{23}F=\begin{cases}
\Z_{p} &\quad p= 691 \\
\Z_{p}\oplus \Z_{p}/(\tau(p)-1) &\quad p\ne 691,
\end{cases}
\]
where $\tau$ is the Ramanujan $\tau$ function, defined by
\[
q\prod_{n=1}^{\infty}(1-q^{n})^{24}=\sum_{k=1}^{\infty} \tau(k)q^{k}.
\]
For $p\ne 691$ the torsion subgroup of $\pi_{23}F$ has order
determined by the $p$-adic valuation of $(\tau(p)-1)$.  The only
primes less that $35,000$ for which $\tau(p)\equiv 1\mod p$ are $11$,
$23$, and $691$.  It is not known whether or not $\tau(p)\equiv 1\mod
p$ holds for infinitely many primes.

The fiber of $\log_{p}^{\tmf}$ seems to store quite a bit of
information about the spectrum of the Atkin operator and $p$-adic
properties of modular forms.  Investigating its homotopy type looks
like interesting prospect for algebraic topology.

\bibliographystyle{amsplain}
\providecommand{\bysame}{\leavevmode\hbox to3em{\hrulefill}\thinspace}
\providecommand{\MR}{\relax\ifhmode\unskip\space\fi MR }

\providecommand{\MRhref}[2]{%
  \href{http://www.ams.org/mathscinet-getitem?mr=#1}{#2}
}
\providecommand{\href}[2]{#2}

\label{lastpage}

\end{document}